\documentclass[11pt]{amsart}

\newtheorem*{theorem}{Theorem}
\newtheorem*{bak}{Anti-calculus Proposition}
\newtheorem*{schwarz}{Schwarz lemma}

\usepackage{url}

\usepackage{tikz}
\usetikzlibrary{patterns}

\usepackage{hyperref}

\renewcommand{\Im}{\mathop{\textrm{Im}}}
\renewcommand{\Re}{\mathop{\textrm{Re}}}

\begin{document}

\title{Julius and Julia:\\
Mastering the Art of the Schwarz Lemma}

\author{Harold P. Boas}
\date{23 December 2009}
\address{Department of Mathematics\\
Texas A\&M University\\
College Station, TX 77843-3368}
\email{boas@math.tamu.edu}
\urladdr{http://www.math.tamu.edu/~boas/}

\begin{abstract}
This article discusses classical versions of the Schwarz lemma at the boundary of the unit disk in the complex plane. The exposition includes commentary on the history, the mathematics, and the applications.
\end{abstract}

\maketitle

\section{Introduction.}
Despite teaching complex analysis for a quarter century,
I still didn't know Jack about the Schwarz lemma. That's what I found out after a graduate student buttonholed me at the door of the printer room.

``Oh, professor, I noticed while studying for the qualifying exam that the book of Bak and Newman attributes the `anti-calculus proposition' to Erd\H{o}s, 
but your father's book cites Loewner. How come?''

Since I was on the verge of sending to the production staff the source files of the new edition of my father's \emph{Invitation to Complex Analysis}~\cite{invitation}, the question demanded my instant attention. I soon located the following statement in 
\emph{Complex Analysis}~\cite{bak} by the late Donald~J. Newman (1930--2007) and his student Joseph Bak.

\begin{bak}[{\cite[p.~78]{bak}}, attributed to Erd\H{o}s]
Suppose \(f\) is analytic throughout a closed disc and assumes its maximum modulus at the boundary point~\(\alpha\). Then \(f'(\alpha)\ne 0\) unless \(f\)~is constant.
\end{bak}

Here the word ``analytic'' (I shall use the classical synonym ``holomorphic'') means that the function~\(f\) is represented locally by its Taylor series about each point. To say that \(f\)~is holomorphic at a boundary point of the domain means that \(f\)~extends to be holomorphic in a neighborhood of that point.

Bak and Newman's title for the proposition refers to the familiar principle from real
calculus that the derivative must vanish at an interior point where a differentiable function attains a maximum. But the name seems inapposite to me,
for the derivative need not vanish at a \emph{boundary} extremum. The attribution too is suspect; perhaps Newman learned the result from Paul Erd\H{o}s (1913--1996), but the statement seems to have been known several years before Erd\H{o}s was born.

I do not mean to criticize the book of Bak and Newman. Indeed, the book is an excellent one according to the following criteria: the glowing review in this \textsc{Monthly} \cite{zalcman}, the existence of a Greek translation, and the observation that the copy in my university's library is falling apart from overuse. 

In any case, the proposition is an arresting statement---though far from optimal---and I learned a lot trying to discover who first proved it, and why. My aim here is to clarify both the history and the mathematics of the proposition and its thematic offshoots by taking a tour through some variations of the Schwarz lemma.  By the end of the tour, you will understand both my title\footnote{Note to readers from future generations: When I wrote this article, movie theaters were showing the motion picture \emph{Julie and Julia}, an exaltation of American cultural icon Julia Child, author of \emph{Mastering the Art of French Cooking}.} 
and the pun in my opening sentence. 

But the story is not entirely lighthearted. We mathematicians are engaged in a social enterprise, and world events can shake the foundations of our ivory tower. I include some biographical remarks about the participants in my tale, for they lived in ``interesting times,''
and the nonmathematical aspects of their lives should be remembered along with their theorems.

\section{How the lemma got its name.}
\label{sec:canonical}
Some preliminary normalizations will expose the key issue in the proposition from the introduction. Replacing \(f(z)\) by \(f(z+c)\), where \(c\)~is the center of the initial disk, one might as well assume from the start that the disk is centered at the origin. Having made that simplification, consider the function~\(F\) that sends a point~\(z\) in the standard closed unit disk \(\{\,z\in\mathbb{C}: |z| \le 1\,\}\) to the value \( f(\alpha z)/f(\alpha)\). This function~\(F\) maps the unit disk into itself; fixes the boundary point~\(1\), where the maximum modulus of~\(F\) is attained; and has the property that
\begin{equation}
\label{eq:rescaling}
F'(1) = \alpha f'(\alpha)/f(\alpha).
\end{equation}
Thus the proposition reduces to a property of the derivative of a holomorphic self-mapping of the unit disk at a boundary fixed point.

Every beginning course in complex analysis contains a version of the following statement about the derivative of a holomorphic self-mapping of the unit disk at an \emph{interior} fixed point.

\begin{schwarz}
Suppose\/ \(f\) is a holomorphic function in the open unit disk,\/ \(f(0)=0\), and the modulus of\/~\(f\) is bounded by\/~\(1\). Then either\/ \(|f'(0)|<1\) and\/ \(|f(z)|<|z|\) when\/ \(0<|z|<1\), or else\/ \(f\)~is a rotation (in which case\/ \(|f'(0)|=1\) and\/ \(|f(z)|=|z|\) for all\/~\(z\)).
\end{schwarz}

This unassuming statement has turned out to be extremely fruitful, spawning numerous generalizations, at least two books \cite{wirths,dineen}, and research that continues to this day. Influential extensions by Ahlfors~\cite{ahlfors} and by Yau~\cite{yau} to complex manifolds satisfying suitable curvature hypotheses are widely cited in the current literature. MathSciNet and Zentralblatt MATH list twenty articles
in the recent five-year period 2001--2005 with 
``Schwarz lemma'' in the title. 

The lemma is named for 
the German mathematician Hermann Amandus Schwarz (1843--1921), who is the eponym of various other analytical objects, including the Cauchy--Schwarz inequality, the Schwarz--Christoffel formula in conformal mapping, the Schwarzian derivative, and the Schwarz reflection principle. Every calculus student ought to know about Schwarz's example~\cite[pp.~309--311]{schwarz} showing that the surface area of a cylinder is not equal to the supremum of the areas of inscribed polygonal surfaces. 

The Greek--German mathematician
Constantin Carath\'eodory (1873--1950) named and codified the Schwarz lemma in \cite[p.~110]{cara1912} (having previously used the lemma in \cite{cara1905} and~\cite{cara1907} without giving it a title). 
The autobiographical notes in Carath\'eodory's collected papers 
include the following account of his rediscovery of the lemma~\cite[p.~407]{cara-V}; the date of the indicated events is 1905, the year after Carath\'eodory defended his dissertation on the calculus of variations.
\begin{quote}
Shortly after the end of the semester, P.~Boutroux came to G\"ottingen for a few days, and during a walk he told me about his attempt to simplify Borel's proof of Picard's theorem, which was all the rage then. Boutroux had noticed that this proof works only because of a remarkable rigidity in conformal mapping, but he was unable to quantify the effect. I could not get Boutroux's observation out of my thoughts, and six weeks later, I was able to prove Landau's sharpened version of Picard's theorem in a few lines by means of the proposition that nowadays is called the Schwarz Lemma.
I derived this lemma by using the Poisson integral; then I learned from Erhard Schmidt, to whom I communicated my discovery, not only that the statement is already in Schwarz but that it can be proved in a very elementary way. Indeed, the proof that Schmidt told me can hardly be improved. This is how I arrived at another line of research besides the calculus of variations.\footnote{The free translation from German into English is mine.}
\end{quote}

Just reading the names in this paragraph sends shivers down my spine: truly, giants walked the earth in those days. 
Carath\'eodory himself subsequently became a leading mathematician in Germany, just behind David Hilbert (1862--1943) in stature. Indeed, at the time of the blow-up in 1928 between Hilbert and L.~E.~J. Brouwer (1881--1966) that forced the reorganization of the editorial board of \emph{Mathematische Annalen} \cite{dalen}, Carath\'eodory was one of the chief editors of the journal, the others being Hilbert himself, 
Albert Einstein (1879--1955), and the managing editor,
Otto Blumenthal (1876--1944).

Thus Carath\'eodory had no need to drop names, but what a list of names is contained in his little story! \'Emile Borel (1871--1956) has a lunar crater named after him and is memorialized by such terms as Borel set, Heine--Borel theorem, and Borel--Cantelli lemma; \'Emile Picard (1856--1941) has over six thousand mathematical descendants listed in the Mathematics Genealogy Project and is remembered for his penetrating applications of Schwarz's method of successive approximations, for Picard's theorems in complex analysis, for the Picard group, and so on; Edmund Landau (1877--1938) has been an idol to generations of analysts and number theorists for his elegant books; Erhard Schmidt (1876--1959), a great friend of Carath\'eodory from their student days together in Berlin (where they attended lectures by Schwarz), has both
Hilbert--Schmidt operators and the Gram--Schmidt orthonormalization process named after him.

And what of the French scholar
Pierre Boutroux? His considerable
scientific contributions are largely overshadowed by the work of two prominent relatives: his father was the distinguished philosopher and historian of science \'Emile Boutroux (1845--1921), while his uncle was the incomparable mathematician 
Henri Poincar\'e (1854--1912). Nonetheless, Pierre Boutroux (1880--1922) was himself a sufficiently prominent mathematician that Henry Burchard Fine (1858--1928) hired him at Princeton University in 1913, although Boutroux soon left to join the French army. He survived World War~I but died young, outliving his father by only a few months.

Now the special case of the lemma that
Schwarz handled appears in the notes \cite[pp.~109--111]{schwarz} from his lectures 
during 1869--70, but as far as I know, these notes did not see print until the publication of Schwarz's collected works in 1890. 
In 1884, Poincar\'e too proved a special case (using an argument similar to the modern one), which he singled out as a fundamental lemma 
in a paper on fuchsian groups \cite[p.~231]{poincare}.
Schmidt's elegant proof of the Schwarz lemma, now canonical, is to apply the maximum principle to the function \(f(z)/z\) after filling in the removable singularity at the origin.

A standard illustration \cite[\S69]{cara-conformal} of the power of the Schwarz lemma is the  following proof of Liouville's theorem stating that \emph{the only bounded entire functions} (functions holomorphic in the entire complex plane) \emph{are constants}. After making affine linear transformations of the variables in both the domain and the range,\footnote{The corresponding argument in~\cite{osserman-notices} elides the translation in the range.} one can assume that the entire function~\(f\) fixes the point~\(0\) and has modulus bounded by~\(1\); the goal now is to show that such a function is identically equal to~\(0\). For an arbitrary positive number~\(R\), the function that takes \(z\) to \(f(Rz)\) satisfies the hypothesis of the Schwarz lemma, so \(|f(Rz)|\le |z|\) when \(|z|<1\). Replace~\(Rz\) by~\(w\) to see that \(|f(w)|\le |w|/R\) when \(|w|<R\). Keeping \(w\)~fixed, let \(R\)~tend to infinity to deduce that \(f(w)=0\) for every point~\(w\) in the complex plane. 

\section{Change of base point.}
\label{sec:base}
The original version of the Schwarz lemma applies to a holomorphic function from the unit disk into itself that fixes the origin. If instead the function fixes some other interior point~\(\zeta\), can something still be said?

The answer is simple to find if one knows that the disk is ``symmetric'' in the sense that there is a self-inverse holomorphic mapping of the disk that interchanges two specified points inside the  disk. For present purposes, it suffices to have an explicit formula when one of the two points is the origin. Let the second point inside the disk be~\(\zeta\), and define a mapping~\(\varphi_{\zeta}\) as follows:
\begin{equation*}
\varphi_{\zeta}(z) = \frac{\zeta-z}{1-\overline{\zeta}z}.
\end{equation*} 
Evidently \(\varphi_{\zeta}\) interchanges \(0\) with~\(\zeta\), and \(\varphi_{\zeta}\) maps the unit disk into itself (and the boundary of the disk to the boundary) because, as a routine calculation shows,
\begin{equation}
\label{mobius}
\left| \frac{\zeta-z}{1-\overline{\zeta}z} \right|^{2} = 1- \frac{ (1-|\zeta|^{2}) (1-|z|^{2})} {\left|1-\overline{\zeta}z\right|^{2}}.
\end{equation}
The composite function \(\varphi_{\zeta}\circ \varphi_{\zeta}\) maps the unit disk into itself, fixing both \(0\) and~\(\zeta\), so the case of equality in the Schwarz lemma implies that \(\varphi_{\zeta}\circ \varphi_{\zeta}\) is the identity function; that is, the function~\(\varphi_{\zeta}\) is its own inverse.

The function \(\varphi_{\zeta}\) is the required symmetry of the disk that swaps the points~\(0\) and~\(\zeta\). To interchange two general points \(\eta\) and~\(\zeta\) of the unit disk, one can use the composite function \(\varphi_{\zeta} \circ \varphi_{\varphi_{\zeta}(\eta)} \circ \varphi_{\zeta}\).

Now if \(f\) maps the unit disk into itself, fixing~\(\zeta\), then \(\varphi_{\zeta}\circ f\circ \varphi_{\zeta}\) maps the unit disk into itself, fixing~\(0\). 
The Schwarz lemma implies that the derivative \( (\varphi_{\zeta}\circ f\circ \varphi_{\zeta})'(0)\) has modulus no greater than~\(1\). By the chain rule, this derivative equals \(\varphi_{\zeta}'(\zeta)f'(\zeta)\varphi_{\zeta}'(0)\). But \(\varphi_{\zeta}'(\zeta)\varphi_{\zeta}'(0)=1\), since \(\varphi_{\zeta}\circ \varphi_{\zeta}\) is the identity function. The conclusion is that \(|f'(\zeta)|\le 1\), regardless of the location of the interior fixed point~\(\zeta\).

More generally, if the point~\(\zeta\) is arbitrary (not necessarily fixed by~\(f\)), then the Schwarz lemma applies to the composite function 
\(\varphi_{f(\zeta)}\circ f\circ \varphi_{\zeta}\), so
\begin{equation*}
|\varphi_{f(\zeta)}'(f(\zeta)) f'(\zeta) \varphi_{\zeta}'(0)| \le 1,
\end{equation*}
and computing the derivatives shows that
\begin{equation}
\label{eq:schwarz-base}
|f'(\zeta)| \le \frac{1-|f(\zeta)|^{2}} {1-|\zeta|^{2}}.
\end{equation}
The other conclusion of the Schwarz lemma is that
\begin{equation*}
|\varphi_{f(\zeta)}\circ f\circ \varphi_{\zeta}(z)| \le |z| \qquad \text{when $|z|<1$},
\end{equation*}
or, since \(\varphi_{\zeta}\) is self-inverse,
\begin{equation}
\label{eq:picky}
|\varphi_{f(\zeta)}\circ f(z)| \le |\varphi_{\zeta}(z)| \qquad \text{when $|z|<1$}.
\end{equation}
Written out explicitly, this inequality says the following:
\begin{equation}
\label{eq:pick}
\left|\frac{f(\zeta)-f(z)}{1-\overline{f(\zeta)}f(z)}\right| \le 
\left|\frac{\zeta-z}{1-\overline{\zeta}z}\right|.
\end{equation}

Carath\'eodory already mentions this idea of composing with a holomorphic automorphism of the disk in the paper where he names the Schwarz lemma 
\cite[\S6]{cara1912}, but the first source I know where inequality~\eqref{eq:pick} is written out
is \cite[p.~71]{julia} by
the French mathematician Gaston Julia (1893--1978).
In Julia's long essay (about which more later), 
the inequality is part of a preliminary section in which he collects various auxiliary results. 

Inequality~\eqref{eq:pick} is commonly called the ``Schwarz--Pick lemma'' after the Austrian mathematician Georg Pick (1859--1942), who wrote an influential paper \cite{pick} on the subject two years before Julia's work. Pick's innovation is not so much the inequality (which he does not write explicitly) as the point of view, decisive for subsequent developments in the field: his interpretation is that holomorphic mappings decrease the hyperbolic, non-Euclidean distance in the unit disk (the so-called Poincar\'e metric). A monotonically increasing function of the quantity on the right-hand side of inequality~\eqref{eq:pick}, the hyperbolic distance is invariant under rotations and under all the mappings~\(\varphi_{\zeta}\).

Pick is remembered today also for Nevanlinna--Pick interpolation; for Pick's theorem about the area of lattice polygons; and for his numerous mathematical descendants, the majority of them through his mathematical grandson Lipman Bers (1914--1993), who finished his dissertation in Prague in 1938 just in time to flee the country before the German occupation of 1939. At that time Pick had long since retired, and he was trapped in Prague; he died soon after being interned in the Theresienstadt concentration camp.

\begin{figure}
  \begin{tikzpicture}[thick]
\draw[dashed] (0,0) circle (3cm);
\fill (0,0) node[anchor=north east]{$0$} circle(2pt);
\fill (3,0) node[anchor=west]{$1$} circle (2pt);
  \begin{scope}
    \pgftransformrotate{30}
\fill (2.25,0) node[anchor=east]{$\zeta$} circle (2pt);
\draw[] (1.9636,0) circle (0.7636);
\draw[] (1.44,0) circle (1.44cm);
\draw[] (2.186,0) circle (0.34cm);
\draw (0.92624,0) circle (2.017cm);
    \end{scope}
\end{tikzpicture}
\caption{Some circles with non-Euclidean center~\(\zeta\).}
\label{fig:circles}
\end{figure}
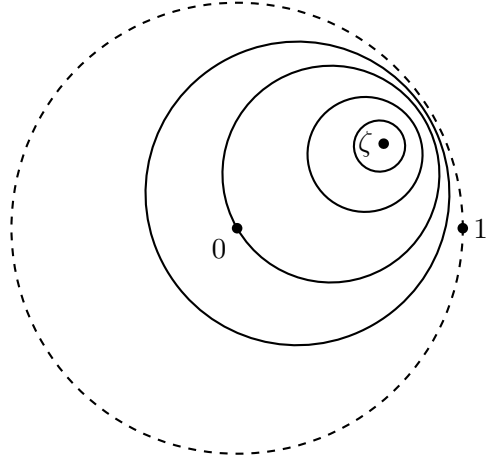

When \(\zeta\)~is a fixed point of~\(f\), Pick's observation provides a simple pictorial representation of
inequality~\eqref{eq:pick}. Setting \(|\varphi_{\zeta}(z)|\) equal to some constant less than~\(1\) restricts the point~\(z\) to lie on a certain Euclidean circle. This circle is also a non-Euclidean circle, and its non-Euclidean center is the point~\(\zeta\). When \(f(\zeta)=\zeta\), inequality~\eqref{eq:pick} says that \(f\) maps points inside this circle to other points inside the circle: that is, \(f\)~contracts the disk bounded by the circle. Figure~\ref{fig:circles} illustrates several such circles that \(f\)~contracts.

\section{Straightening out the problem.}
\label{sec:straight}
There is nothing essential about working on the standard unit disk; one can easily transform the preceding results to apply to a disk with an arbitrary center and an arbitrary radius. It may be less obvious that one can equivalently work on the upper half-plane (thought of as a disk with a boundary point at
infinity); it will be useful to know that one can switch back and forth between the two settings.

One learns in a first course in complex analysis that if \(\zeta\)~is a point in the upper half-plane (that is, \(\Im\zeta >0\)), and 
\begin{equation*}
\psi_{\zeta}(z) = \frac{z-\zeta}{z-\overline{\zeta}},
\end{equation*} 
then \(\psi_{\zeta}\) is an invertible holomorphic function of~\(z\) mapping the upper half-plane onto the unit disk, taking \(\zeta\) to~\(0\). In fact, 
\begin{equation}
\label{eq:UHPidentity}
\left| \frac{z-\zeta}{z-\overline{\zeta}}\right|^{2} = 1- \frac{ 2(\Im z)(\Im\zeta)} {\left|z-\overline{\zeta}\right|^{2}},
\end{equation}
so \(\psi_{\zeta}\) does map the upper half-plane into the unit disk; one can see that the map is bijective by exhibiting the inverse:
\begin{equation*}
\psi_{\zeta}^{-1}(w) = \frac{\zeta-\overline{\zeta}w} {1-w}, \qquad
\text{and} \qquad
\Im \frac{\zeta-\overline{\zeta}w} {1-w} = \frac{(\Im \zeta) (1-|w|^{2})} { |1-w|^{2}}.
\end{equation*}

Now if \(f\) is a holomorphic function that maps the upper half-plane into itself, and \(\zeta\)~is an arbitrary point in the upper half-plane, then the composite function \(\psi_{f(\zeta)} \circ f \circ \psi_{\zeta}^{-1}\) maps the unit disk to itself, fixing the point~\(0\). The derivative \((\psi_{f(\zeta)} \circ f \circ \psi_{\zeta}^{-1})'(0)\) equals \(\psi_{f(\zeta)}'(f(\zeta)) f'(\zeta)/\psi_{\zeta}'(\zeta)\), and one computes that \(\psi_{\zeta}'(\zeta) = (2i\Im\zeta)^{-1}\), so the Schwarz lemma implies that
\begin{equation}
\label{eq:schwarz-upper}
|f'(\zeta)| \le \frac{\Im f(\zeta)}{\Im \zeta} \qquad
\text{when $\Im \zeta >0$.}
\end{equation}
Absolute-value signs are not needed on the right-hand side, for both the numerator and the denominator are positive real numbers. If \(\zeta\)~happens to be a fixed point of~\(f\), then \(|f'(\zeta)|\le 1\), just as in the case of a disk of finite radius. 

The bound \eqref{eq:schwarz-upper} on \(|f'(\zeta)|\) looks a bit simpler than the corresponding bound~\eqref{eq:schwarz-base} in the case of the unit disk; an analogous simplification occurs for the inequality involving the function~\(f\) itself. In the setting of the preceding paragraph, the Schwarz lemma implies additionally that 
\begin{equation}
\label{eq:half-plane}
|\psi_{f(\zeta)}\circ f(z)| \le |\psi_{\zeta}(z)|,
\qquad \text{or} \qquad
\left| \frac{f(z)-f(\zeta)} {f(z)-\overline{f(\zeta)}}\right|  
\le \left| \frac{z-\zeta}{z-\overline{\zeta}}\right|
\end{equation}
when \(\zeta\) and~\(z\) lie in the upper half-plane.
As far as I know, Julia was the first to write 
this half-plane version of the Schwarz--Pick lemma explicitly \cite[p.~71]{julia}; an equivalent form, in view of equation~\eqref{eq:UHPidentity}, is that
\begin{equation}
\label{eq:UHPpick}
\frac{ (\Im f(z))(\Im f(\zeta))} {\left|f(z)-\overline{f(\zeta)}\right|^{2}} \ge
\frac{ (\Im z)(\Im\zeta)} {\left|z-\overline{\zeta}\right|^{2}}.
\end{equation}

Raised in Algeria, Julia went to Paris on a scholarship, but World War~I interrupted his education. Serving on the front lines in 1915, he was grievously wounded, shot in the middle of his face. 
Julia finished his dissertation during the long recuperation from his injury and the subsequent unsuccessful operations to repair the damage (for the rest of his life, he wore a leather strap across his face to hide the missing nose~\cite{mendes}). Then he studied iteration of rational functions in the plane and wrote a monograph on the subject for a competition announced by the Acad\'emie des Sciences; this work won the 1918 Grand Prix. Julia's prize-winning essay subsequently fell into relative obscurity until it was revived in the 1970s by Benoit Mandelbrot in his celebrated work on fractals. Nowadays computer-generated images of ``Julia sets'' are ubiquitous.

Of primary interest here is 
the section in Julia's article devoted to the Schwarz lemma \cite[pp.~67--78]{julia}, which is only a small part of the nearly 200-page work. Julia perhaps recognized the independent significance of this part of his treatise, for he later published a version of this section as a stand-alone note~\cite{julia-bis}. 

\section{Fixed point on the boundary.}
\label{sec:boundary}
The discussion so far is, for Julia, child's play. But now he is
ready to present his extension of the Schwarz lemma. He supposes that a nonconstant holomorphic function~\(f\) maps a disk into itself and has a fixed point \emph{on the boundary}.
For this situation to make sense, the function~\(f\) ought to be defined on the boundary, or at least at the fixed point. Julia assumes initially that \(f\)~is holomorphic on the closed disk, and he finds it convenient to work on the half-plane model of a disk, as in the preceding section. After making a translation, one might as well take the boundary fixed point to be the origin. 

Now \(f(z)=f'(0)z + O(z^{2})\). Evidently \(f'(0)\) cannot be equal to~\(0\), for if \(f(z)\) were to behave locally near the origin like a square (or a higher power), then \(f\)~would multiply angles by a factor of~\(2\) (or more), contradicting the hypothesis that \(f\)~maps the upper half-plane into itself. Multiplication by the nonzero complex number~\(f'(0)\) corresponds to a dilation by \(|f'(0)|\) composed with a rotation by \(\arg f'(0)\). Since \(f\)~maps the upper half-plane into itself, the rotation must be the identity. Thus \(f'(0)\)~is actually a positive real number. By the method of the preceding section, this conclusion carries over to the setting of the standard unit disk.

The preceding observations prove the proposition from the introduction, and more. In view of the chain-rule computation~\eqref{eq:rescaling}, not only does a nonconstant holomorphic function~\(f\) on a disk centered at the origin have nonzero derivative at a boundary point~\(\alpha\) where the modulus~\(|f|\) takes a maximum, but in addition the expression \(\alpha f'(\alpha)/f(\alpha)\) is a positive real number. 

The quantity \(\alpha f'(\alpha)/f(\alpha)\) has a compelling geometric interpretation. The function~\(f\) maps the boundary circle onto a curve that lies inside a circle centered at the origin with radius \( |f(\alpha)|\) and touches this circle at the point~\(f(\alpha)\). The complex number \( i\alpha f'(\alpha)\) represents a vector tangent to the curve at \(f(\alpha)\), and an outward normal vector is \(-i\) times this quantity, or \(\alpha f'(\alpha)\). But the curve is tangent to the circle at~\(f(\alpha)\), so the normal vector to the curve points in the same direction as the radius vector of the circle, which is \( f(\alpha)\). Thus the positivity of the ratio \( \alpha f'(\alpha)/f(\alpha)\) expresses the tangency of the image curve with the enclosing circle.

Although Julia carries out most of his analysis in 
the upper half-plane, he states explicitly that the positivity of the derivative at a boundary fixed point holds for an arbitrary disk.
Actually, the whole development so far is as easy as pie for Julia, but now comes the pi\`ece de r\'esistance.

The Schwarz--Pick lemma for the upper half-plane yields the inequality~\eqref{eq:UHPpick} for values of the function~\(f\); what more can be said if the boundary point~\(0\) is a fixed point? Simply taking the limit in~\eqref{eq:UHPpick} as \(\zeta\) approaches~\(0\) yields a trivial inequality, for both sides reduce to~\(0\). But useful information can be extracted from~\eqref{eq:UHPpick} by considering the \emph{rate} at which the two sides approach~\(0\). 

Keeping \(z\) fixed, replace \(f(\zeta)\) by \(f'(0)\zeta + O(\zeta^{2})\), divide both sides of~\eqref{eq:UHPpick} by \(\Im \zeta\), and let \(\zeta\) approach~\(0\) along the imaginary axis. What results is the following fundamental inequality of Julia:
\begin{equation}
\label{eq:julia}
f'(0)\frac{\Im f(z)}{|f(z)|^{2}} \ge \frac{\Im z} {|z|^{2}} 
\qquad \text{when $\Im z>0$}.
\end{equation}
Julia actually writes the equivalent statement that
\begin{equation}
\label{eq:julia2}
\Im \frac{1}{f(z)} \le \Im \frac{1}{f'(0)z}.
\end{equation}
Since  \(\Im (1/z) = - \Im(z)/|z|^{2}\), with a minus sign,
inequality~\eqref{eq:julia2} has the direction reversed relative to inequality~\eqref{eq:julia}. Julia's form of the inequality reveals that if equality holds at even one point in the upper half-plane, then the harmonic function 
\begin{equation*}
\Im\left( \frac{1}{f(z)} -  \frac{1}{f'(0)z} \right)
\end{equation*}
attains a maximum and so is constant; therefore \(f\)~is a M\"obius transformation (a rational function of degree~\(1\)).

Notice that the preceding derivation of Julia's inequality does not really require the function~\(f\) to be defined on the boundary. As long as there is a sequence of points~\(\zeta_n\) in the upper half-plane tending to~\(0\) such that \(f(\zeta_n)\to 0\) and the fraction \( \Im (f(\zeta_n))/\Im(\zeta_n)\) stays bounded above, the same argument implies that inequality~\eqref{eq:julia} holds with \(f'(0)\) replaced by \( \liminf_{n\to\infty} \{ \Im (f(\zeta_n))/\Im(\zeta_n) \}\).

\begin{figure}[t]
\begin{tikzpicture}[xscale=2,yscale=2,thick]
\small
\draw (0,0.83333) circle (0.40825cm);
\draw (2,2) circle (0.979796cm);
\draw (-2,0) -- (3.5,0);
\fill (0,0.5) circle (1pt);
\path node[above left=-2pt] at (0,0.5) {$\zeta$};
\fill (2,1.2) node[above left=-2pt]{$f(\zeta)$} circle (1pt);
\filldraw (0,0.83333) circle(1pt) -- (0,1.24158) node[pos=0.5,anchor=west]{$r$};
\filldraw (2,2) circle (1pt) -- (2,2.979796) node[pos=0.5,anchor=west]{$R$};
\draw[dashed] (-1.42857,0) -- (2,1.2);
\fill (-1.42857,0) node[anchor=north]{$p$} circle (1pt);
\draw[dashed] (0,0) -- (0,0.5) node[pos=0.5,anchor=west]{$\Im \zeta$};
\draw[dashed] (2,0) -- (2,1.2) node[pos=0.5,anchor=west]{$\Im f(\zeta)$};
\draw[dashed] (-1.42857,0) -- (2,2);
\draw[dashed] (-1.42857,0) -- (2,2.979796);
\end{tikzpicture}
\caption{Two congruent non-Euclidean disks.}
\label{fig:julia}
\end{figure}
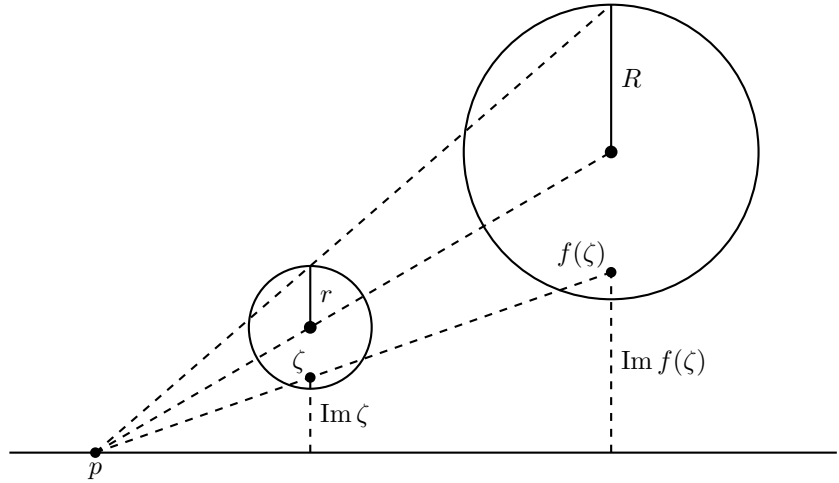

\begin{figure}[bht]
\begin{tikzpicture}[xscale=2,yscale=2,thick]
\small
\begin{scope}
\pgfsetfillpattern{crosshatch dots}{blue!40} 
\filldraw (0,0) arc (-90:270:0.979796cm) arc (270:-90: 0.40825cm);
\end{scope}
\begin{scope}
\filldraw[fill=red!20] (0,0) arc (-90:270:0.40825cm);
\end{scope}
\filldraw (0,0.40825) circle (1pt) -- (0,0.8165) node[pos=0.5,anchor=west]{$r$};
\filldraw (0,0.979796) circle (1pt) -- (0,1.959572) node[pos=0.5,anchor=west]{$R$};
\draw (-1,0) -- (1,0);
\begin{scope}[yshift=0.40825cm]
\fill (120:0.40825cm) node[anchor=south east]{$z$} circle (1pt);
\end{scope}
\begin{scope}[yshift=1cm]
\filldraw (0.4,0) node[anchor=south west]{$f(z)$} circle (1pt);
\end{scope}
\fill (0,0) node[anchor=north]{$0$} circle (1pt);
\draw (2,1) node{\normalsize $\xrightarrow[w\mapsto -1/w]{\text{via inversion}}$};
\begin{scope}[xshift=4cm]
\begin{scope}
\fill[fill=red!20] (-1, 1.225) -- (1, 1.225) -- (1,2) -- (-1,2) -- (-1, 1.225);
\end{scope}
\begin{scope}
\pgfsetfillpattern{crosshatch dots}{blue!40} 
\fill (-1, 0.5103) -- (1, 0.5103) -- (1, 1.225) -- (-1, 1.225) -- (-1, 0.5103);
\end{scope}
\draw (-1,0) -- (1,0);
\draw (-1, 1.225) -- (1, 1.225) node[anchor=west]{$\dfrac{1}{2r}$};
\draw (-1, 0.5103) -- (1, 0.5103) node[anchor=west]{$\dfrac{1}{2R}$};
\fill (0,0) node[anchor=north]{$0$} circle (1pt);
\fill (0.32824, 1.225) node[anchor=south]{$-1/z$} circle (1pt);
\fill (-0.3448,0.86207) node[anchor=north]{$-1/f(z)$} circle (1pt);
\end{scope}
\end{tikzpicture}
\caption{Inverting the limiting configuration.}
\label{fig:limit}
\end{figure}
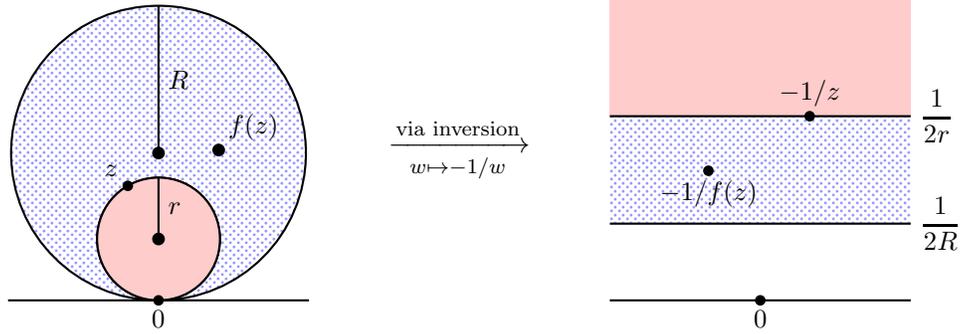

The proof that I have given is not Julia's. His illuminating geometric argument is the following. Consider in the upper half-plane two disks with equal non-Euclidean radius, one centered at a point~\(\zeta\) and the other centered at~\(f(\zeta)\); see Figure~\ref{fig:julia}.
The Schwarz--Pick lemma implies that the image of the first disk under the
holomorphic function~\(f\) lies inside the second disk.
A suitable dilation with center at the indicated point~\(p\) on the real axis takes the point~\(\zeta\) to the point~\(f(\zeta)\). Being an invertible holomorphic mapping, this dilation preserves the non-Euclidean distance and so maps the first disk precisely onto the second disk. Accordingly, the ratio \(R/r\) of the Euclidean radii equals the dilation factor \( (\Im f(\zeta))/(\Im \zeta)\). Now suppose \(\zeta\)~lies on the imaginary axis, and let \(\zeta\)~move vertically down to the origin while the Euclidean radius~\(r\) stays fixed. The variable quantity \(R/r\) then tends to a limit that equals the limit of \( (\Im f(\zeta))/(\Im \zeta)\), or~\(f'(0)\). See Figure~\ref{fig:limit}, where the point~\(\zeta\) has disappeared in the limit, and I have introduced a new, arbitrary point~\(z\).
After mapping each point to the negative of its reciprocal, one sees in the right-hand part of the figure that dilating the point \(-1/f(z)\) by the factor \(R/r\) gives a point that lies above the horizontal line with height equal to \(\Im(-1/z)\). But \(R/r=f'(0)\), so \(f'(0)\Im(-1/f(z)) \ge \Im(-1/z)\), which is Julia's inequality again.

Julia has no need to write out 
the equivalent inequality for a holomorphic function~\(f\) mapping the unit disk into itself, fixing the point~\(1\), but
a routine computation (conjugating with the mapping~\(-\psi_{i}\) in the notation of Section~\ref{sec:straight}) leads to the following result:
\begin{equation}
\label{eq:julia-wolff}
f'(1)\ge \frac{1-|z|^{2}}{|1-z|^{2}} \biggm/ \frac{1-|f(z)|^{2}} {|1-f(z)|^{2}} \qquad \text{when $|z|<1$}.
\end{equation}
Carath\'eodory makes this inequality explicit in~\cite{cara1929}.

For more on Carath\'eodory's life, consult the meticulously researched biography \cite{georgiadou}, but watch out for the author's occasional lapses in English. For instance, the author says that Carath\'eodory married ``his aunt'' \cite[p.~51]{georgiadou}, when actually the relationship was that of second cousins once removed. (Carath\'eodory's great-grandfather Antonios was the brother of Carath\'eodory's wife's grandfather Stephanos.)

\section{Multiple fixed points.}
\label{sec:twofer}
The case of equality in the Schwarz lemma implies that the only holomorphic self-mapping of the unit disk with two interior fixed points is the identity function. On the other hand, there are nontrivial holomorphic self-mappings of disks that fix an interior point together with one or more
boundary points. For example, the mapping that sends~\(z\) to~\(z^{3}\) takes the unit disk into itself and fixes both the interior point~\(0\) and the boundary points \(1\) and~\(-1\). Transferring such a situation to the upper half-plane and applying Julia's fundamental inequality~\eqref{eq:julia} with \(z\)~equal to the interior fixed point
shows immediately that the derivative at each boundary fixed point is greater than or equal to~\(1\); in fact, the derivative is strictly greater than~\(1\) unless the function is the identity.

In the case of a disk centered at the origin, the rescaling trick of Section~\ref{sec:canonical} leaves the origin fixed, so the chain-rule calculation~\eqref{eq:rescaling} implies the following statement.

\begin{theorem}[after Julia]
If a holomorphic function~\(f\) (not identically~\(0\)) maps a closed disk centered at~\(0\) into itself, fixing~\(0\), and if the modulus~\(|f|\) attains a maximum on the closed disk at the boundary point~\(\alpha\), then \(\alpha f'(\alpha)/f(\alpha)\) is a real number greater than or equal to~\(1\).
\end{theorem}

That this theorem follows immediately from Julia's inequality is not a new observation \cite[p.~524]{ruscheweyh}, but I have not determined to my satisfaction where the theorem was first written down in this form. 

The unscaled case when \(\alpha=1=f(\alpha)\)
appears as problem~291 in the famous problem book of P\'olya and Szeg\H{o} \cite{polyaszego}, the first edition of which was published in German in 1925. In the solution, the authors indicate that the problem was contributed\footnote{In the author index, the entry that points to the solution has the page number in \emph{italics}, and a note at the beginning of the index indicates that this typographical convention indicates ``original contributions.''} by Charles Loewner.\footnote{They actually say ``K.~L\"owner,'' his name at the time. Loewner was caught in Prague when the Germans occupied Czechoslovakia in 1939, but nonetheless he succeeded in emigrating.
His student Lipman Bers explains \cite[p.~vii]{loewner-papers} that Loewner's ``decision to Americanize a name already well known to mathematicians the world over was a sign of the 46-year-old scholar's determination to begin a new life in the United States.''} 
It is no surprise that Loewner was thinking along these lines, for he was a student of Pick, and Loewner addresses related questions in the article \cite{loewner} in which he introduces his famous differential equation for studying schlicht functions. 

W.~K. Hayman, writing after
World War~II, makes the related observation that 
``the points where \(f(z)\) attains its maximum modulus lie among the points where \(z_{0} f'(z_{0})/f(z_{0})\) is real''
\cite[p.~137]{hayman}, and he attributes the idea (though not the formula) to Blumenthal's remarkable 1907 paper~\cite{blumenthal}, a study of the differentiability properties of the maximum-modulus function. I can extract at least the proposition in the introduction from that paper \cite[p.~223, formula~(7)]{blumenthal}.

I learned from the notes in R.~B. Burckel's encyclopedic treatise \cite[p.~216]{burckel} that the scaled formulation of the theorem, essentially as I have given it above, appears in a 1971 paper by I.~S. Jack \cite[Lemma~1]{jack}; the paper is the unique item listed in MathSciNet for this author. The subsequent literature on univalent functions commonly refers to the theorem as ``Jack's lemma.''\footnote{Some authors use the name ``Clunie--Jack lemma'' on the strength of Jack's comment, ``we use a method due to Professor Clunie.'' It seems that Clunie himself accepted this designation while also acknowledging antecedents \cite[p.~139]{clunie}.} 
On the basis of this contribution, I.~S. Jack has an impressive ratio of citations to published articles (more than 250 hits in Google Scholar).

The theorem can be sharpened by invoking Schmidt's trick. If \(f\) maps the unit disk into itself, fixing both \(0\) and~\(1\), then introduce the function~\(g\) such that \(g(z)=f(z)/z\) when \(z\ne 0\) and \(g(0)=f'(0)\). For a suitable choice of angle~\(\theta\), the function \(e^{i\theta} \varphi_{g(0)}\circ g\) again maps the unit disk into itself, fixing both \(0\) and~\(1\), so the derivative at~\(1\) is no smaller than~\(1\). Translating this conclusion to a statement about~\(f\) by a routine calculation shows that
\begin{equation*}
f'(1) \ge \frac{2}{1+|f'(0)|}.
\end{equation*}
Since \(|f'(0)|\le 1\), this inequality refines the assertion that \(f'(1)\ge 1\).
This improvement was obtained in \cite{unkelbach} by Helmut Unkelbach (1910--1968); it was rediscovered more than half a century later in~\cite{osserman}.

\section{Julius has a point.}
\label{sec:wolff}
A pleasing counterpart to the preceding section occurs when all the fixed points lie on the boundary. In this case there is a distinguished fixed point at which the derivative is a positive number less than or equal to~\(1\); at all the other fixed points, the derivative is larger than~\(1\) and moreover at least as large as the reciprocal of the derivative at the special fixed point.

This proposition is a corollary of results of the Dutch mathematician\footnote{This Julius Wolff should not be confused with the German orthopedic surgeon of the same name who lived 1836--1902, nor with the German novelist and lyricist who lived 1834--1910.}
Julius Wolff (1882--1945). 
Here is a concrete example:
\begin{equation*}
f(z)=
\left( \frac{z-\frac{2}{3}} {1- \frac{2}{3}z} \right)^{3}.
\end{equation*}
This function~\(f\) maps the unit disk into itself, fixing the four boundary points \(-1\), \(+1\), and \((9\pm 5i\sqrt{7}\,)/16\); computing shows that
\begin{equation*}
f'(-1)=3/5, \qquad f'(1)=15, \qquad \text{and} \qquad f'((9\pm 5i\sqrt{7}\,)/16) = 12/5.
\end{equation*}

The easy part of the proposition is the statement 
that the product of derivatives at two distinct boundary fixed points must be at least~\(1\). After transferring the problem to the upper half-plane, one can assume that the two fixed points are \(0\) and~\(1\). Now invoke Julia's inequality~\eqref{eq:julia}, multiplying both sides by the positive quantity \(|f(z)|^{2}/\Im z\) and taking the limit as \(z\)~approaches~\(1\) along a vertical line. 
Since \(f(1)\) (which equals~\(1\)) and \(f'(1)\) are real numbers, 
\begin{equation*}
\frac{\Im f(z)}{\Im z} = \frac{\Im [f(z)-f(1)]}{\Im (z-1)} = \frac{f'(1)\Im(z-1)+O((z-1)^{2})} {\Im(z-1)},
\end{equation*}
so it follows that \(f'(0)f'(1)\ge 1\).
I will show at the end of the proof that the derivatives \(f'(0)\) and \(f'(1)\) cannot both be equal to~\(1\). 

One of the objectives of Wolff's theory is to eliminate the assumption that the function is defined on the boundary. To accommodate this generalized setting, one can declare a boundary point~\(\alpha\) to be a fixed point of~\(f\) when \(\alpha\)~is equal to the limit of \(f(z)\) as \(z\) approaches~\(\alpha\) along the direction normal to the boundary.

Suppose, then, that \(f\)~is a holomorphic function mapping the open upper half-plane into itself. Wolff observes \cite{wolff-mono} that the quotient \( (\Im f(z))/(\Im z)\) increases when \(z\) moves down a vertical line (where \(\Re z\) is constant). Indeed, let~\(z\) be the point in Figure~\ref{fig:julia} at the top of the circle with Euclidean radius~\(r\). 
As observed in Section~\ref{sec:boundary}, 
the image point~\(f(z)\) lies inside or on the indicated circle with Euclidean radius~\(R\). Hence the ratio \((\Im f(z))/(\Im z)\) is less than or equal to the dilation ratio \(( \Im f(\zeta))/(\Im \zeta)\). Alternatively, one can verify Wolff's observation by a calculus argument. If \(y\)~denotes the imaginary part of~\(z\), then
\begin{equation*}
\frac{\partial}{\partial y} \left( \frac{\Im f(z)}{\Im z} \right) = 
\frac{1}{\Im z}\left( \Re f'(z) - \frac{\Im f(z)}{\Im z}\right).
\end{equation*}
But \(\Re f'(z) \le |f'(z)|\), so inequality~\eqref{eq:schwarz-upper} implies that the derivative on the left-hand side is negative or zero. Consequently, when \(\Im z\) decreases to~\(0\) with \(\Re z\) fixed, the ratio \( (\Im f(z))/(\Im z)\) increases either to a positive limit or to~\(+\infty\).

Next suppose that there are no interior fixed points. For this part of the argument, it is convenient to work on the unit disk. 
Wolff \cite{wolff-horo} adapts an elegant 
idea that Arnaud Denjoy (1884--1974) used \cite{denjoy} to simplify the proof of an earlier result of Wolff. When \(f\)~is a holomorphic function mapping the open unit disk into itself, and \(n\)~is a positive integer, consider the function~\(f_{n}\) defined by \(f_{n}(z) = (1-\tfrac{1}{n})f(z)\). There are various ways to see that \(f_{n}\) has a fixed point~\(\alpha_{n}\) inside the unit disk: Brouwer's fixed-point theorem (apply it to a closed disk of radius close to~\(1\)), the contraction-mapping theorem (the non-Euclidean disk is a complete metric space), or Rouch\'e's theorem (the function \(f_{n}(z)-z\) has the same number of zeroes inside the disk as \(z\)~does: namely, one). 
Since \(\alpha_{n} = f_{n}(\alpha_{n})\), no subsequence of the~\(\alpha_{n}\) can converge to an interior point~\(\alpha\) of the disk, else \(\alpha\) would equal \(f(\alpha)\), contrary to the hypothesis that \(f\)~has no interior fixed point. Thus \(|\alpha_{n}|\to 1\) as \(n\to \infty\).

\begin{figure}
\begin{tikzpicture}[thick]
\draw[dashed] (0,0) circle (3cm);
\fill (0,0) node[anchor=north east]{$0$} circle(2pt);
\fill (3,0) node[anchor=west]{$1$} circle (2pt);
  \begin{scope}
    \pgftransformrotate{30}
    \draw (3,0) arc (0:360:2cm);
\draw (3,0) arc (0:360:1.5cm);
  \draw (3,0) arc (0:360:1cm);
  \draw (3,0) arc (0:360:0.5cm);
  \fill (3,0) node[anchor=south west]{$\alpha$} circle (2pt);
    \end{scope}
\end{tikzpicture}
\caption{Some horocycles.}
\label{fig:horo}
\end{figure}
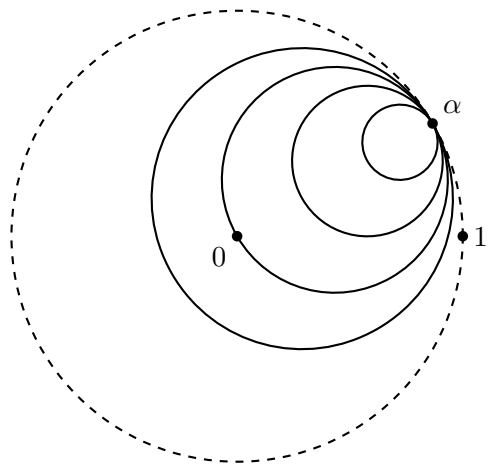

Version~\eqref{eq:picky} of the Schwarz--Pick lemma says that
\begin{equation*}
|\varphi_{\alpha_n}( f_n(z))| \le |\varphi_{\alpha_n}(z)| \qquad \text{when $|z|<1$}.
\end{equation*}
As observed at the end of Section~\ref{sec:base}, 
what this inequality means geometrically is that for each~\(z\), the point \(f_n(z)\) lies inside or on the circle~\(C_n(z)\) that passes through~\(z\) and has non-Euclidean center~\(\alpha_n\). Some subsequence
of the~\(\alpha_n\) converges
to a boundary point~\(\alpha\)
(after the fact, one can deduce that the whole sequence converges), and passing to the limit leads to the conclusion that the point \(f(z)\) lies inside or on the circle~\(C(z)\) passing through~\(z\) and tangent to the unit circle at~\(\alpha\). In other words, the function~\(f\) contracts circles that are tangent to the unit circle at~\(\alpha\).
Figure~\ref{fig:horo} illustrates some of these circles, called horocycles. By identity~\eqref{mobius}, 
a point~\(w\) lies on the circle \(C_{n}(z)\) if and only if
\begin{equation*}
\frac{1-|w|^{2}} {|1-\overline{\alpha_{n}}w|^{2}} = \frac{1-|z|^{2}} {|1-\overline{\alpha_{n}}z|^{2}},
\end{equation*}
so \(w\) lies on the horocycle \(C(z)\) if and only if
\begin{equation*}
\frac{1-|w|^{2}} {|1-\overline{\alpha}w|^{2}} = \frac{1-|z|^{2}} {|1-\overline{\alpha}z|^{2}}.
\end{equation*}

The special boundary point~\(\alpha\) is known as the Denjoy--Wolff point.
The linking of the names seems appropriate: Wolff had two notes in the \emph{Comptes Rendus} on January~4 and January~18 of 1926, Denjoy's follow-up was on January~25, and Wolff revisited the subject on April~12 and September~13. 
Although French, Denjoy held a professorship for five years at Utrecht University, and Wolff succeeded him in that position.

Now it is convenient to switch back to the equivalent setting of the upper half-plane, placing the Denjoy--Wolff point at~\(0\). Since \(f\)~contracts horocycles tangent to the real axis at~\(0\), the limit of \(f(z)\) as \(z\)~moves down the imaginary axis to~\(0\) equals~\(0\); that is, the origin is a boundary fixed point in the generalized sense. Moreover, \(\Im f(z) \le \Im z\) when \(z\)~is on the imaginary axis. Hence the limit of the ratio \( ( \Im f(z))/(\Im z)\) at~\(0\) along the imaginary axis (shown to exist at the beginning of the discussion) is a positive number~\(\beta\) less than or equal to~\(1\). If \(f\)~extends to be continuously differentiable at~\(0\), then \(\beta=f'(0)\). In any case, the remark preceding inequality~\eqref{eq:julia-wolff} shows that Julia's inequality~\eqref{eq:julia} holds with~\(\beta\) in the role of~\(f'(0)\).

One technical point needs attention to complete the proof of the proposition that I stated at the beginning of this section. Could there be two fixed points on the real axis, say \(0\) and~\(1\), such that both of the derivatives \(f'(0)\) and \(f'(1)\) (more generally, the corresponding limiting ratios~\(\beta\))
are equal to~\(1\)? The answer is no, because Julia's argument from Section~\ref{sec:boundary} implies that the holomorphic function~\(f\) would contract horocycles based at~\(1\) as well as horocycles based at~\(0\). Two such horocycles 
having Euclidean radius~\(1/2\) touch at one point, and this point would be an interior fixed point of~\(f\), contradicting the assumption that all the fixed points lie on the boundary.

Wolff actually shows considerably more than I have stated. In particular, he proves in~\cite{wolff-mono} that the difference quotient \( (f(z)-f(\alpha))/(z-\alpha)\)
has a limit when \(z\)~approaches the distinguished boundary point~\(\alpha\) within an arbitrary triangle contained in the disk and with a vertex at~\(\alpha\). This property was soon rediscovered by 
Landau and Georges Valiron (1884--1955) in~\cite{landau-valiron} and by Carath\'eodory, who systematized the theory in~\cite{cara1929} and introduced the terminology ``angular derivative'' (Winkelderivierte) for arbitrary boundary fixed points (not just the Denjoy--Wolff point).

In summary, a holomorphic self-mapping of a disk (not equal to the identity function) has a distinguished fixed point---either in the interior or on the boundary---at which the derivative has modulus less than or equal to~\(1\). The (angular) derivative exceeds~\(1\) at all the other fixed points.

Wolff returned to the topic in 1940, proving in \cite{wolff-product} a refined statement about the product of angular derivatives at two boundary fixed points. That year the German army invaded the Netherlands, and Wolff was forced out of the university six months later. He continued publishing nonetheless and collaborated with Blumenthal, who had taken refuge in the Netherlands in 1939 after six years of persecution in Germany. 
(All four of Blumenthal's grandparents were Jewish.)
But in April of 1943, the order came to deport to the concentration camps all remaining persons in Utrecht of Jewish heritage. 
Blumenthal perished in November 1944 in Theresienstadt~\cite{felsch}, and Wolff three months later in Bergen-Belsen~\cite{digital}.

\section{Boundary uniqueness.}
The case of equality in the Schwarz lemma can be viewed as saying that if a holomorphic mapping from a disk into itself looks like the identity function to first order at an interior point, then the mapping \emph{is} the identity function. The analogous statement at a boundary point fails: the function that sends~\(z\) to \(z-\tfrac{1}{4}(z-1)^3\) maps the unit disk into itself (as a straightforward calculus calculation shows), looks like the identity
at the boundary point~\(1\) not only to first order but even to second order, but is not the identity.

Apparently nobody thought of looking for a general positive result about uniqueness at the boundary during the first half of the twentieth century. It turns out that third-order tangency to the identity suffices. The following statement comes from a paper mainly devoted to multidimensional boundary uniqueness problems.

\begin{theorem}[Burns and Krantz {\cite[Theorem 2.1, p.~662]{burnskrantz}}]
Let \(\phi\) be a holomorphic function mapping the unit disk into itself such that
\begin{equation*}
\phi(z)=z + O((z-1)^{4}) \qquad \text{as $z\to 1$}.
\end{equation*}
Then \(\phi\) is the identity function.
\end{theorem}

The proof in \cite{burnskrantz} is short but sophisticated; the authors invoke both Hopf's lemma and the Herglotz representation of a positive harmonic function as the Poisson integral of a positive measure. Incidentally, the proposition in the introduction is a very special case of Hopf's lemma, which itself is a small part of an influential paper~\cite{hopf} 
by the Austrian mathematician Eberhard Hopf (1902--1983) on the maximum principle for subsolutions of elliptic partial differential equations.

Remarkably, Julia's fundamental inequality yields a short and elementary proof of the uniqueness theorem. 
The argument is most conveniently written in the equivalent setting of a holomorphic function~\(f\) mapping the upper half-plane into itself such that \(f(z)=z+O(z^{4})\) as \(z\to 0\). Let \(g(z)\) denote \( z^{-1} - f(z)^{-1}\). Then \(g(z) = O(z^{2})\) as \(z\to 0\), and Julia's inequality~\eqref{eq:julia2} implies that \(\Im g(z) \ge 0\) for all~\(z\) in the upper half-plane. It is easy to see that these two conditions force \(g\) to be the \(0\)~function.

If \(g\) happens to extend holomorphically to a neighborhood of the point~\(0\), then the argument is the same as in Section~\ref{sec:boundary}: if \(g\)~is not identically~\(0\), then \(g\)~is multiplying angles at the origin by a factor of~\(2\) (or more), so \(g\) maps some points of the upper half-plane into the lower half-plane, contradicting that \(\Im g(z)\ge 0\).

When \(g\)~lives only in the open upper half-plane (the general case), the conclusion follows directly from Harnack's inequality, a standard tool in potential theory, 
which says in particular that a \emph{positive} harmonic function in a disk cannot tend to~\(0\) at the boundary along a radius at a faster rate than the distance to the boundary. Since \(\Im g(z)\) tends to~\(0\) at least quadratically  as \(z\)~goes down the imaginary axis toward the origin, Harnack's inequality implies that the harmonic function \(\Im g\), although nonnegative, cannot be strictly positive. Thus \(\Im g\) is a harmonic function that assumes its extreme value~\(0\), so the maximum principle (applied to the function~\(-\Im g\)) shows that \(\Im g(z)\) is identically equal to~\(0\). Hence \(g(z)\)~too is identically equal to~\(0\), so \(f(z)\) is identically equal to~\(z\).

Harnack's inequality appears in a book \cite[\S19, p.~62]{harnack} by the Baltic German mathematician Axel Harnack (1851--1888). Harnack's statement is that if \(u\)~is a harmonic function of fixed sign in a disk of radius~\(r\), and \(u_{0}\) is the value of~\(u\) at the center point, then the value of~\(u\) at a point at distance~\(\rho\) from the center is between \(u_{0}(r+\rho)/(r-\rho)\) and \(u_{0}(r-\rho)/(r+\rho)\). This inequality follows easily from the Poisson integral representation for harmonic functions.

\section{Conclusion.}
In this article I have discussed the classical Schwarz lemma, the boundary versions discovered by Gaston Julia and Julius Wolff, and some applications. This circle of ideas is known as Julia--Wolff theory or Julia--Wolff--Carath\'eodory theory. Readers who wish to go further with the mathematics might start with the expository papers \cite{burckel-iterating} and~\cite{osserman-notices}, Carath\'eodory's book for a treatment of angular derivatives \cite[\S\S298--304]{cara-functions}, the book~\cite{abate} concerning implications for iteration theory both in one dimension and in higher dimensions, and the book~\cite{dineen} for applications in infinite-dimensional analysis. Some entry points to the historical literature are 
the books \cite{segal} and \cite{fleeing} for the lives of mathematicians in a desperate time;
the book~\cite{audin} for an in-depth treatment of the 1918 Grand Prix, its antecedents, and its aftermath; and the book~\cite{alexander} for the early development of iteration theory.

Remember that graduate student who asked me the question? He passed the qualifying exam with flying colors. 

\providecommand{\bysame}{\leavevmode\hbox to3em{\hrulefill}\thinspace}

\end{document}